\documentclass{svproc}

\usepackage{graphicx}
\usepackage{subcaption}
\usepackage{amsmath}
\usepackage{bm} 

\graphicspath{{figs}}

\begin{document}

\mainmatter             

\title{Operator-difference approximations on two-dimensional merged Voronoi-Delaunay grids}

\titlerunning{Approximations on two-dimensional merged Voronoi-Delaunay grids}  

\author{Petr N. Vabishchevich}

\authorrunning{Petr N. Vabishchevich}

\tocauthor{Petr N. Vabishchevich}

\institute{
	Lomonosov Moscow State University, Moscow 119234, Russia,\\
	\email{vab@cs.msu.ru}
}

\maketitle             

\begin{abstract}
Formulating boundary value problems for multidimensional partial derivative equations in terms of invariant operators of vector (tensor) analysis is convenient. Computational algorithms for approximate solutions are based on constructing grid analogs of vector analysis operators.
This is most easily done by dividing the computational domain into rectangular cells when the grid nodes coincide with the cell vertices or are the cell centers.
Grid operators of vector analysis for irregular regions are constructed using Delaunay triangulations or Voronoi partitions.
This paper uses two-dimensional merged Voronoi-Delaunay grids to represent the grid cells as orthodiagonal quadrilaterals.
Consistent approximations of the gradient, divergence, and rotor operators are proposed.
On their basis, operator-difference approximations for typical stationary scalar and vector problems are constructed.
\keywords{finite volume method, vector analysis operators, Voronoi polygons, Delaunay triangulation}
\end{abstract}

\section{Introduction}\label{sec:1}

In difference methods \cite{Samarskii1989,Strikwerda2004}, an approximate solution to boundary value problems is sought by approximating the solution at the nodes of the computational grid.
On irregular computational grids, difference approximations are constructed using the balance method (integral-interpolation method)  \cite{Samarskii1989,TikhonovSamarskii1961}.
In this case, the grid problem is an integral consequence of the mathematical model equations (conservation law) for separate parts of the computational domain.
This approach is known as the finite volume method (see, e.g., \cite{FVM,Generalized_difference_methods}).

It is convenient to formulate boundary value problems for multidimensional partial derivative equations using vector (tensor) analysis operators that are invariant regarding the choice of the coordinate system.
Computational algorithms for approximate solving boundary value problems on arbitrary computational grids are constructed using discrete analogs of such operators.
The grid operators of vector analysis should inherit the main properties of differential operators related to fulfilling the fundamental integral relations.
Such consistent approximations are the basis of MD (Mimetic Discretization) technology \cite{shashkov2018conservative,da2014mimetic}.

Consistent approximations of vector analysis operators are relatively easy to construct on rectangular grids \cite{castillo2013mimetic}.
When arbitrary grids are used, success is achieved by applying the method of support operators \cite{SamarskiiKoldobaPoveschenkoFavorskii1996,shashkov2018conservative}.
The grid analogs of vector analysis operators on irregular computational grids are most easily obtained using Delaunay triangulations and Voronoi partitions \cite{vabishchevich2005finite}.

When partitioning the computational domain into cells, the standard approach involves choosing cell vertices as nodes of the computational grid or cell centers.
A two-grid technique for difference approximations using these two computational grids simultaneously is proposed in  \cite{frjazinov1975,samarskiiUMN}.
The grid problem is considered the direct sum of the corresponding Hilbert spaces of the grid functions defined on the two grids.
For scalar and vector boundary value problems, it is possible to construct operator-difference approximations on a common unified rectangular grid \cite{vabishchevich2024gvm}.
In the present paper, we consider the problems of vector analysis operator approximations on a two-dimensional merged Voronoi-Delaunay grid.
In this case, the computational domain is partitioned into orthodiagonal quadrilaterals.
Consistent approximations of gradient, divergence, and rotor operators are constructed.
The use of grid operators in vector analysis is illustrated in typical scalar and vector boundary value problems.

\section{Grids and grid functions}\label{sec:2}

Among unstructured grids, we can distinguish Delaunay triangulation, for which the Voronoi partitioning is the dual.
We define a general MVD (Merged Voronoi-Delaunay) grid whose nodes are both the Delaunay triangulation nodes and the Voronoi partition nodes.

\begin{figure}[htbp]
\centering
\includegraphics[width=0.75\linewidth]{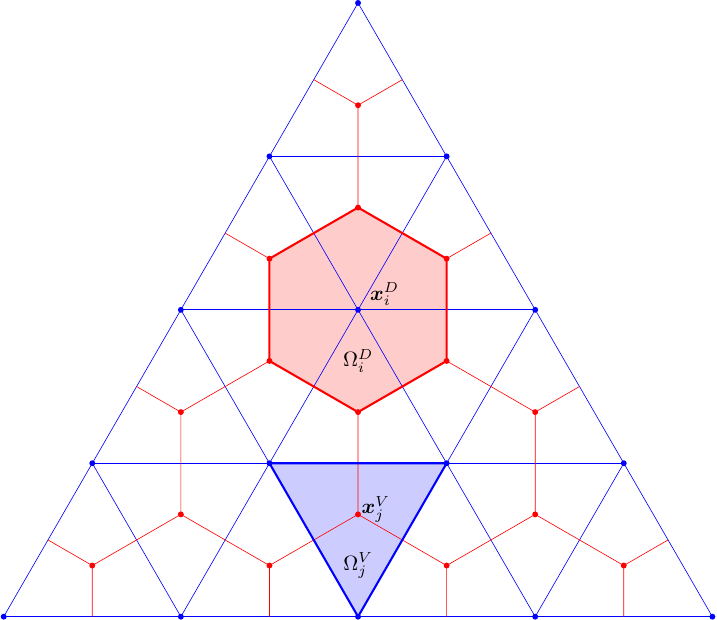} 
\caption{Nodes and grid cells in Delaunay triangulation and Voronoi partitioning}
\label{f-1}
\end{figure}

We assume the computational domain is a convex two-dimensional polygon $\Omega$ with boundary $\partial \Omega$.
In the domain $\overline{\Omega} = \Omega \cup \partial \Omega$ we have nodes $\bm x_i^D, \ i = 1,2, \ldots, M_D$ of the $D$-grid
\[
\omega^D = \{ \bm x \ | \ \bm x = \bm x_i^D, \  \bm x_i^D \in \overline{\Omega}, \ i = 1,2, \ldots, M_D \} .
\]
The set of nodes on the boundary of $\partial \Omega$ includes the vertices of the polygon $\Omega$.

We associate each node $\bm x_i^D, \ i = 1,2, \ldots, M_D$ with a part of the computational domain $V_i^D$, which is a Voronoi polygon or a part of it belonging to $\Omega$.
The Voronoi polygon for an individual node $i$ is the set of points lying closer to this node $i$ than to all others $i^\pm= 1,2, \ldots, M_D$ and therefore
\[
\Omega_i^D = \{ \bm x \ | \ \bm x \in \Omega , \ |\bm x - \bm x_i^D | < |\bm x - \bm x_{i^\pm}^D |, \ i^\pm= 1,2, \ldots, M_D \}, \quad i = 1,2, \ldots, M_D ,
\]
where $| \cdot |$ is the Euclidean distance.
Let $\bm x_j^V, \ j = 1,2, \ldots, M_V$ are vertices of $V_i^D, \ i = 1,2, \ldots, M_D$.
We will assume that all Voronoi polygons vertices lie inside the computational domain $\Omega$ or on its boundary $\partial \Omega$.
The grid $V$ is defined by the set of nodes (vertices of the Voronoi diagram polygon) $\bm x_j^V, \ j = 1,2, \ldots, M_V$:
\[
\omega^V = \{ \bm x \ | \ \bm x = \bm x_j^V, \ \ \bm x_j^V \in \overline{\Omega}, \ j = 1,2, \ldots, M_V \} .
\]
The cells of the $D$-grid are the parts of the computational domain $\Omega_i^D, \ i = 1,2, \ldots, M_D$:
\[
\overline{\Omega} = \bigcup_{i=1}^{M_D} \overline{\Omega}_i^D, \ \overline{\Omega}_i^D = \Omega_i^D \cup \partial \Omega_i^D,
\Omega_{i^\pm}^D \cap \Omega_i^D = \emptyset , \ i^\pm\neq i, \ i^\pm, \ i = 1,2, \ldots, M_D .
\]

Each vertex $\bm x_j^V, \ j = 1,2, \ldots, M_V$  of the Voronoi polygon is associated with the triangle constructed by the appropriate nodes contacting the Voronoi polygons.
The cells of the $V$-grid are the parts of the computational domain 
\[
 \Omega_j^V = \{ \bm x \ | \ \bm x \in \Omega , \ |\bm x - \bm x_j^V | < |\bm x - \bm x_{j^\pm}^V |, \ j^\pm= 1,2, \ldots, M_V \}, \quad j = 1,2, \ldots, M_V ,
\] 
so that 
\[
 \overline{\Omega} = \bigcup_{j=1}^{M_V} \overline{\Omega}_j^V, \ \overline{\Omega}_j^V = \Omega_j^V \cup \partial \Omega_j^V,
 \Omega_{j^\pm}^V \cap \Omega_j^V = \emptyset , \ j^\pm\neq j, \ j^\pm, j = 1,2, \ldots, M_V .
\] 

These triangles determine the Delaunay triangulation --- a dual triangulation to the Voronoi diagram. 
We assume that the Delaunay triangulation and the Voronoi partitioning are sufficiently regular.
In particular, each node of the $\bm x_j^V$ $V$-grid lies inside $\overline{\Omega}_j^V , \ j = 1,2, \ldots, M_V$.
The introduced notations are illustrated by Fig.~\ref{f-1} for the nodes and cells of the two-dimensional computational grid.  

\begin{figure}[htbp]
\centering
\includegraphics[width=0.75\linewidth]{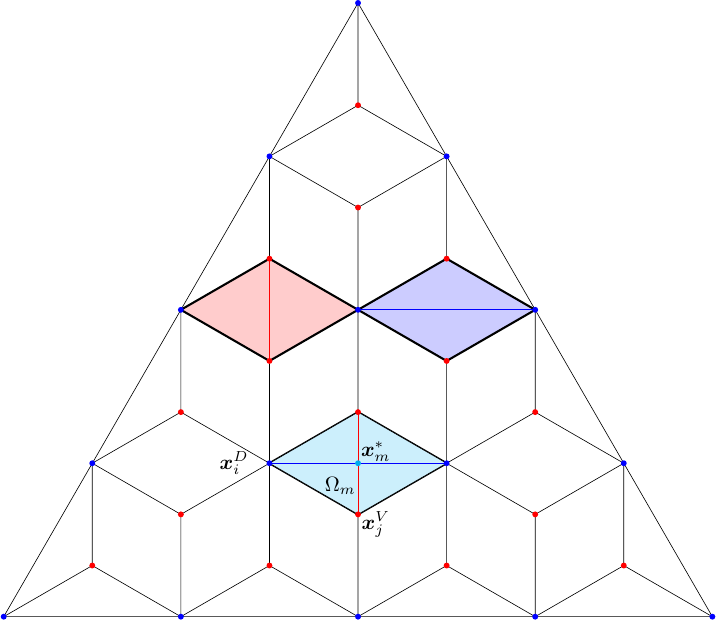} 
\caption{MVD grid nodes and MVD grid cells in a two-dimensional domain}.
\label{f-2}
\end{figure}

We consider a general grid with both $D$-grid nodes and $V$-grid nodes.
In the two-dimensional case we consider, the cells of the MVD grid are the orthodiagonal quadrilaterals $\Omega_m, \ m = 1,2, \ldots, M$ (triangles adjoining the boundary, see Fig.~\ref{f-2}).
In the $\Omega$ computational domain, a grid is used
\[
\overline{\Omega} = \bigcup_{m=1}^{M} \overline{\Omega}_m, \ \overline{\Omega}_m = \Omega_m \cup \partial \Omega_m,
\Omega_{m^\pm} \cap \Omega_m = \emptyset , \ m^\pm \neq m, \ m^\pm , m = 1,2, \ldots, M .
\]

For the functions $u(\bm x), \ \bm x \in \Omega$, we define a Hilbert space $\mathcal{H} = L_2(\Omega)$ with scalar product and norm
\[
(u,v) = \int_{\Omega}u(\bm x) v(\bm x) \, d \bm x,
\quad \|u\| = (u,u)^{1/2} .
\]
The Hilbert spaces of the grid functions on the introduced grids are similarly introduced.
The grid functions are defined at the grid nodes with which a part of the computational domain (control volume) is associated.

For the grid functions defined at the nodes of the Delaunay triangulation $\bm x_i^D, \ i = 1,2, \ldots, M_D$, we define a finite-dimensional Hilbert space $H(\omega^D)$, in which
\[
(y,v)_D = \sum_{\bm x \in \omega^D} y(\bm x) \, v(\bm x) S^D(\bm x) ,
\quad \|y\|_D = (y,y)_D^{1/2} .
\]
The control volumes for the Delaunay triangulation are Voronoi polygon so that
\[
S^D(\bm x) = \operatorname{meas} (\Omega_{i}^D) ,
\quad \bm x = \bm x_{i}^D \in \omega^D .
\]

\begin{figure}[htbp]
\centering
\includegraphics[width=0.75\linewidth]{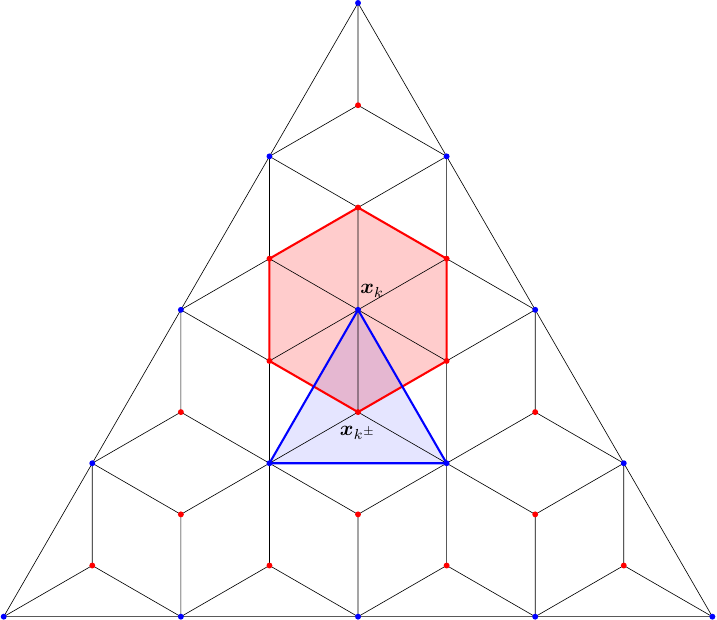} 
\caption{MVD grid nodes and control volumes }
\label{f-3}
\end{figure}

We similarly define a finite-dimensional Hilbert space $H(\omega^V)$ for the grid functions described at the vertices of Voronoi polygons $\bm x_j^V, \ j = 1,2, \ldots, M_V$:
\[
(y,v)_V = \sum_{\bm x \in \omega^V} y(\bm x) \, v(\bm x) S^V(\bm x) ,
\quad \|y\|_V = (y,y)_V^{1/2} .
\]
For the nodes of the $\omega^V$ grid, the control volumes are the cells of the Delaunay triangulation:
\[
S^V(\bm x) = \operatorname{meas} (\Omega_{j}^V) ,
\quad \bm x = \bm x_{j}^V \in \omega^V .
\]

We have $\omega = \omega^D \cup \omega^V$ for the merged grid.
Let us define the space $H(\omega)$ with scalar product and norm
\[
(y,v) = \sum_{\bm x \in \omega} y(\bm x) \, v(\bm x) S(\bm x) ,
\quad \|y\| = (y,y)^{1/2}
\]
for the grid functions $y(\bm x), \ \bm x \in \omega$ on the MVD grid.
The previously introduced control volumes for the nodes of the $\omega^D$ grid and the $\omega^V$ grid cover the computational domain twice (see Fig.~\ref{f-3}) and therefore for the nodes $\bm x_k, \ k = 0,1, \ldots, M_{VD}$ of the $\omega$ grid we obtain
\[
 S(\bm x) = \frac{1}{2} \left \{ 
 \begin{array}{ll}
 S^D(\bm x) , & \bm x \in \omega^D , \\
 S^V(\bm x) , & \bm x \in \omega^V.
 \end{array}
 \right . 
\] 

For the introduced MVD grid, grid functions are defined at the nodes of cells --- node-center approximation.
Note also the possibility of using the second class of grids when the grid nodes are the centers of cells --- cell-center approximation.
The cell center is chosen as the intersection of the orthogonal diagonals of the quadrilateral cell $\Omega_m, \ m = 1,2, \ldots, M$ (Fig.~\ref{f-3}).
We denote the cell-center nodes of the grid by $\bm x_m^*, \ m = 1,2, \ldots, M$, and the grid itself by $\omega^*$.
For the grid functions defined at nodes $\bm x_m^*, \ m = 1,2, \ldots, M$, we define a finite-dimensional Hilbert space $H(\omega^*)$ in which
\[
(y,v)_* = \sum_{\bm x \in \omega^*} y(\bm x) \, v(\bm x) S^*(\bm x) ,
\quad \|y\|_* = (y,y)_*^{1/2} .
\]
For the control volumes, we have
\[
S^*(\bm x) = \operatorname{meas} (\Omega_{m}) ,
\quad \bm x = \bm x_{m}^* \in \omega^* .
\]
In our case, the main grid is $\omega$, whose nodes are the vertices of cells $\Omega_m, \ m = 1,2, \ldots$, and the auxiliary grid is the grid
$\omega^*$ with nodes inside the cells.

\section{Grid operators of vector analysis}

We construct grid analogs of invariant vector analysis operators (independent of the choice of the coordinate system).
The key feature of our approximations is related to using the main $\omega$ for scalar quantities and the auxiliary grid $\omega^*$ --- for vector quantities.

Using the VAGO (Vector Analysis Grid Operator) \cite{vabishchevich2005finite} technique, we construct discrete analogs of the vector analysis differential operators: gradient, divergence, and rotor.
In the MD (Mimetic Discretization) approach \cite{shashkov2018conservative,da2014mimetic}, we keep track of the fundamental properties of discrete vector analysis operators.

We will construct the grid operator of the gradient at nodes $\bm x_m^*, \ m = 1,2, \ldots, M$ of the auxiliary grid $\omega^*$ by the grid function given at nodes $\bm x_k, \ k = 0,1, \ldots, M_{VD}$ of the main grid $\omega$.
The equality $\bm w = \operatorname{grad} u$ is matched by
\begin{equation}\label{1}
\bm v (\bm x^*) = \operatorname{grad}_h y(\bm x),
\quad \bm x \in \omega ,
\quad \bm x^* \in \omega^* ,
\end{equation}
where the grid operator of the gradient $\operatorname{grad}_h: H(\omega) \longrightarrow \bm H(\omega^*)$.
For vector grid functions from a finite-dimensional Hilbert space $\bm H(\omega^*)$, let us put
\[
 (\bm y,\bm v)_* = \sum_{\bm x \in \omega^*} \bm y(\bm x) \, \bm v(\bm x) S^*(\bm x) ,
 \quad \|y\|_* = (y,y)_*^{1/2} .
\] 

\begin{figure}[htbp]
\centering
\includegraphics[width=1\linewidth]{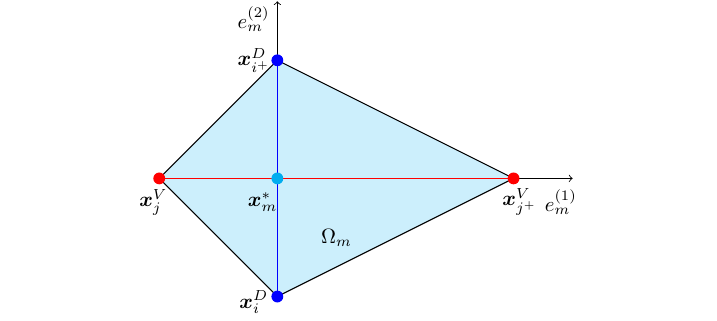} 
\caption{Gradient approximation on MVD grid }
\label{f-4}
\end{figure}

It is convenient to introduce for an individual cell $\Omega_m$ of the MVD grid a local Cartesian coordinate system (see Fig.~\ref{f-4}) centered at node $\bm x_m^*$ and coordinates $e^{(1)}_{m}, e^{(2)}_{m}, \ m = 1,2, \ldots, M$.
We will approximate the individual components of the vector $\bm v (\bm x_m^*) = \{v^{(1)}(\bm x_m^*), v^{(2)}(\bm x_m^*)\}$:
\begin{equation}\label{2}
(\operatorname{grad}_h y) (\bm x_m^*) = \left \{ \frac{y(\bm x_{j^+}^V) - y(\bm x_{j}^V)}{|\bm x_{j^+}^V - \bm x_{j}^V|} ,
\frac{y(\bm x_{i^+}^D) - y(\bm x_{i}^D)}{|\bm x_{i^+}^D - \bm x_{i}^D|} \right \} .
\end{equation}

The approximations (\ref{1}), (\ref{2}) can be obtained by starting from the integral definition of the gradient.
For the control volume $\Omega_{m}$ we have
\begin{equation}\label{3}
\int_{\Omega_{m}}\operatorname{grad} u (\bm x) \, d\bm x = \int_{\partial \Omega_{m}} u(\bm x) \, \bm n \, d \bm x ,
\end{equation}
where $\bm n$ is the external normal to $\partial \Omega_{m}$.
Given this, the grid operator of the gradient at node $\bm x_m^* \in \omega^*$ is defined according to
\begin{equation}\label{4}
(\operatorname{grad}_h y)(\bm x_m^*) = \frac{1}{S^*(\bm x)} \int_{\partial \Omega_{m}} \widetilde{y} (\bm x) \, \bm n \, d \bm x ,
\end{equation}
where $\widetilde{y}(\bm x)$ is a piecewise linear interpolation on parts of the boundary of the control volume $\partial \Omega_{m}$ over the values of the grid function $y(\bm x), \ \bm x \in \omega$.
By integrating over the contour using the quadrature trapezoidal formula on the right-hand side (\ref{4}), we arrive at the approximation (\ref{2}).

Similarly, the grid operator of divergence is constructed based on the integral equality for the corresponding control volume.
Given a vector field $\bm v (\bm x)$ at nodes $\bm x \in \omega^*$, we approximate $y(\bm x) = \operatorname{div}_h \bm v$
at nodes $\bm x \in \omega$.

\begin{figure}[htbp]
\centering
\includegraphics[width=0.7\linewidth]{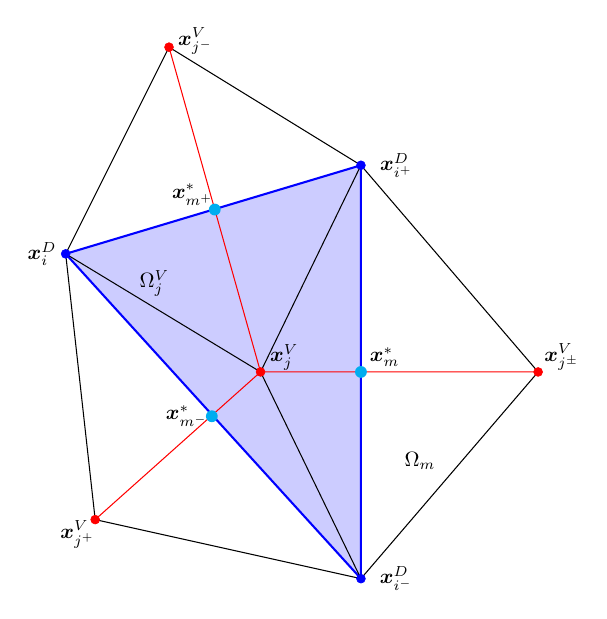} 
\caption{Divergence approximation on the MVD grid}
\label{f-5}
\end{figure}

Let us consider the case when a node of the two-dimensional grid $\omega$ is a Voronoi partition node and a Delaunay triangulation cell acts as the control volume (Fig.~\ref{f-5}).
The grid operator of divergence is constructed based on the integral equality for the corresponding control volume.

Applying the divergence theorem to the control volume $\Omega_j^V$ for the node $\bm x_j^V \in \omega$ gives
\begin{equation}\label{5}
\int_{\Omega_j^V}\operatorname{div} \bm w (\bm x) \, d \bm x = \int_{\partial \Omega_j^V} \bm w \cdot \bm n \, d \bm x .
\end{equation}
Given (\ref{5}), a grid operator divergence $\operatorname{div}_h: \bm H(\omega^*) \longrightarrow H(\omega)$ is defined according to
\begin{equation}\label{6}
(\operatorname{div}_h \bm v) (\bm x_j^V)  = \frac{1}{S(\bm x_j^V)} \int_{\partial \Omega_j^V} \widetilde{\bm v} \cdot \bm n \, d \bm x ,
\end{equation}
where $\widetilde{\bm v}(\bm x)$ is a piecewise constant interpolation of the grid vector function $\bm v(\bm x), \ \bm x \in \omega^*$ on parts of the boundary of the control volume $\partial \Omega_j^V$.

From (\ref{6}), using the notations of Fig.~\ref{f-5}, we obtain the following representation of the grid operator of divergence
\begin{equation}\label{7}
\begin{split}
(\operatorname{div}_h \bm v) (\bm x_j^V) = \frac{1}{S(\bm x_j^V)}  \Big ( & v_n(\bm x_m^*) |\bm x_{i^+}^D - \bm x_{i^-}^D| \\
+ & v_n(\bm x_{m^+}^*) |\bm x_{i^+}^D - \bm x_{i}^D| + v_n(\bm x_{m^-}^*) |\bm x_{i}^D - \bm x_{i^-}^D| \Big ) .
\end{split}
\end{equation}
The divergence approximation is performed using the normal components $v_n (\bm x) = (\bm v(\bm x) \cdot \bm n)$ of the vector $\bm v(\bm x)$ at the nodes of the auxiliary grid $\bm x \in \omega^*$.
Similarly (\ref{7}), we define the grid operator of divergence at the nodes of the Delaunay triangulation.
For example, when computing $ (\operatorname{div}_h \bm v) (\bm x_i^D)$, the expressions used are
\[
\frac{1}{S(\bm x_i^D)} v_n(x_{m^-}^*) |\bm x_{j^+}^V - \bm x_{j}^V|,
\quad \frac{1}{S(\bm x_i^D)} v_n(x_{m^+}^*) |\bm x_{j^-}^V - \bm x_{j}^V|,
\]
where now the normal components $v_n(x_{m^-}^*), \ v_n(x_{m^+}^*)$ of the vector $\bm v(\bm x)$ are defined along the normal to the boundary of the control volume $\Omega_i^D$.

The rotor operator for two-dimensional problems has its peculiarities.
In this case, the vector $\bm u$ in the general case has three components ($\bm u = \{u_1, u_2, u_3\}$), which for the problems we consider depend only on the coordinates $x_1, x_2$: $u_\alpha = u_\alpha(x_1, x_2), \ \alpha = 1,2,3$.
Let's extract two components into the vector $\bm u$:
\[
\bm u = \bm u^{\parallel} + \bm u^{\perp} ,
\quad \bm u^{\parallel} = u_1 \bm e^{(1)} + u_2 \bm e^{(2)},
\quad \bm u^{\perp} = u_3 \bm e^{(3)} ,
\]
where $\bm e^{(\alpha )}, \ \alpha =1,2,3$  are the orthants of the Cartesian coordinate system.
For the vector rotor we have
\begin{equation}\label{8}
\operatorname{rot} \bm u = \operatorname{rot2 \textbf{D}} (\bm u^{\parallel}) \bm e^{(3)} + \operatorname{rot2D} (u^{\perp}).
\end{equation}
Given that
\[
\operatorname{rot} \bm u =  \frac{\partial u_3}{\partial x_2}  \bm e^{(1)} - \frac{\partial u_3}{\partial x_1}  \bm e^{(2)}
+ \left ( \frac{\partial u_2}{\partial x_1} - \frac{\partial u_1}{\partial x_2} \right ) \bm e^{(3)} ,
\]
let us define the rotor operators on the right-hand side (\ref{8}).
For vector $\bm v = \{v_1, v_2\}$ and scalar $w$ we have
\begin{equation}\label{9}
\operatorname{rot2 \textbf{D}} \bm v = \frac{\partial v_2}{\partial x_1} - \frac{\partial v_1}{\partial x_2} ,
\quad \operatorname{rot2D} (w) = \frac{\partial w}{\partial x_2}  \bm e^{(1)} - \frac{\partial w}{\partial x_1}  \bm e^{(2)} .
\end{equation}
We will consider the approximation of $\operatorname{rot2 \textbf{D}}$ and $\operatorname{rot2D}$ separately.

The vector field $\bm u = u \bm e^{(3)}$ will be set in the nodes of the main grid.
With the equality $\bm w = \operatorname{rot} \bm u$ is associated with the approximation
\begin{equation}\label{10}
\bm v (\bm x^*) = \operatorname{rot2D}_h y(\bm x) \bm e^{(3)},
\quad \bm x \in \omega ,
\quad \bm x^* \in \omega^* ,
\end{equation}
so that $\operatorname{rot2D}_h: H(\omega) \rightarrow \bm H(\omega^*)$.
Similarly to the gradient approximation (see Fig.~\ref{4}), if we use a local coordinate system for the individual components of the vector, we have
$\bm v (\bm x_m^*) = \{v^{(1)}(\bm x_m^*), v^{(2)}(\bm x_m^*)\}$:
\begin{equation}\label{11}
(\operatorname{rot2D}_h y) (\bm x_m^*) = \left \{ \frac{y(\bm x_{i^+}^D) - y(\bm x_{i}^D)}{|\bm x_{i^+}^D - \bm x_{i}^D|} ,
- \frac{y(\bm x_{j^+}^V) - y(\bm x_{j}^V)}{|\bm x_{j^+}^V - \bm x_{j}^V|}  \right \} .
\end{equation}
The approximation (\ref{11}) agrees with the approximation (\ref{2}) for the gradient operator.

The vector field in the plane $\bm u = \{u_1, u_2\}$ is defined on an auxiliary grid. 
Given (\ref{9}) and $\operatorname{rot} \bm u = w \bm e^{(3)}$, let's set 
\begin{equation}\label{12}
 y (\bm x) = \operatorname{rot2\textbf{D}}_h \bm v(\bm x^*) ,
 \quad \bm x^* \in \omega^* ,
 \quad \bm x \in \omega.
\end{equation} 
Grid operator $\operatorname{rot2\textbf{D}}_h: \bm H(\omega^*) \rightarrow H(\omega)$ will be constructed based on the integral definition of the rotor.
For example, for the control volume $\Omega_m$ (see Fig.~\ref{5}):
\begin{equation}\label{13}
 \int_{\Omega_{m}} (\operatorname{rot} \bm u \cdot \bm e^{(3)} ) \, d \bm x = 
 \oint_{\partial \Omega_{m}} ( \bm u \cdot \bm \tau ) \, d \bm x , 
\end{equation} 
where $\bm \tau$ is the tangent vector and the contour $\partial \Omega_{m}$ is oriented via the right-hand rule.
The approximation (\ref{13}) leads us to the following expression for the operator $\operatorname{rot2\textbf{D}}_h$:
\begin{equation}\label{14}
\begin{split}
 (\operatorname{rot2\textbf{D}}_h \bm v) (\bm x_j^V) = \frac{1}{S(\bm x_j^V)}  \Big ( & v_\tau(\bm x_m^*) |\bm x_{i^+}^D - \bm x_{i^-}^D| \\\
 + & v_\tau(\bm x_{m^+}^*) |\bm x_{i^+}^D - \bm x_{i}^D| + v_\tau(\bm x_{m^-}^*) |\bm x_{i}^D - \bm x_{i^-}^D| \Big ) 
\end{split}
\end{equation} 
when $\bm v (\bm x) = \{v_n, v_\tau\}$ for $\bm x \in \omega^*$.
The grid approximation (\ref{14}) is compared to the approximation (\ref{7}) for the divergence operator.

\section{Boundary value problems}\label{sec:4} 

The constructed approximations of vector analysis operators on the MVD grid can be used in approximate solutions of boundary value problems.
We give examples of the construction of operator-difference approximations for typical two-dimensional boundary value problems.

In $\Omega$, we are to find the solution to the Dirichlet problem for the second-order elliptic equation
\begin{equation}\label{15}
- \operatorname{div} \big (k(\bm x) \operatorname{grad} u \big ) + c(\bm x) u = f(\bm x),
\quad \bm x \in \Omega ,
\end{equation}
\begin{equation}\label{16}
u(\bm x) = 0,
\quad \bm x \in \partial \Omega ,
\end{equation}
with sufficiently smooth coefficients and right-hand side at $k(\bm x) \geq k_0 = \operatorname{const} > 0$, $c(\bm x) \geq c_0 = \operatorname{const} > 0$.
We define an approximate solution $y(\bm x)$ to the problem (\ref{15}), (\ref{16}) on the grid $\omega$, wherein
\begin{equation}\label{17}
y(\bm x) = 0,
\quad \bm x \in \partial \omega ,
\end{equation}
where $\partial \omega = \omega \cap \partial \Omega$.
The operator-difference approximation of Eq. (\ref{15}) gives
\begin{equation}\label{18}
- \operatorname{div}_h \big (k(\bm x) \operatorname{grad}_h y \big ) + c(\bm x) y = f(\bm x),
\quad \bm x \in \omega - \partial \omega .
\end{equation}
Here the grid operators of gradient and divergence are defined according to (\ref{2}) and (\ref{7}).
On the set of grid functions (\ref{17}), they are consistent in the sense of fulfilling the equality
\[
(\operatorname{grad}_h y, \bm v) + (y, \operatorname{div}_h \bm v) = 0.
\]
Taking this into account, let us write the discrete problem (\ref{17}), (\ref{18}) as an operator equation
\begin{equation}\label{19}
A y = \varphi , \quad A = A^* > 0 ,
\end{equation}
in which the operator $A: H(\omega) \rightarrow H(\omega)$ is
\[
 A y = - \operatorname{div}_h \big (k(\bm x) \operatorname{grad}_h y \big ) + c(\bm x) y.
\]  

The second example of a boundary value problem is related to the determination of $\bm u(\bm x)$ from the conditions of
\begin{equation}\label{20}
 - \operatorname{grad} \big (k(\bm x) \operatorname{div} \bm u \big ) + c(\bm x) \bm u = \bm f(\bm x), 
 \quad \bm x \in \Omega , 
\end{equation} 
\begin{equation}\label{21}
 (\bm u \cdot \bm n) (\bm x) = 0,
 \quad \bm x \in \partial \Omega .
\end{equation} 
The approximate solution $\bm y$ of the problem (\ref{20}), (\ref{21}) is sought at the nodes of the auxiliary grid $\omega^*$.
On the set of grid functions 
\[
 (\bm y \cdot \bm n) (\bm x) = 0, 
 \quad \bm x \in \partial \omega^* , 
\] 
where $\partial \omega^* = \omega^* \cap \partial \Omega$, let's define the operator $A: \bm H(\omega^*) \rightarrow \bm H(\omega^*)$ such that 
\[
 A \bm y = - \operatorname{grad}_h \big (k(\bm x) \operatorname{div}_h \bm y \big ) + c(\bm x) \bm y.
\]  
The boundary value problem (\ref{20}), (\ref{21}) is compared to Eq.
\begin{equation}\label{22}
 A \bm y = \bm \varphi , \quad A = A^* > 0 . 
\end{equation} 

Another example of vector problems involve determining $\bm u(\bm x)$ from conditions
\begin{equation}\label{23}
 \operatorname{rot} \big (k(\bm x)) \operatorname{rot} \bm u \big ) + c(\bm x) \bm u = \bm f(\bm x), 
 \quad \bm x \in \Omega , 
\end{equation} 
\begin{equation}\label{24}
 (\bm u \times \bm n) (\bm x) = 0,
 \quad \bm x \in \partial \Omega .
\end{equation} 
Among the two-dimensional problems (\ref{23}), (\ref{24}), we distinguish separately the case when $\bm u(\bm x) = u(\bm x) \bm e^{(3)}$.
For the scalar function $u(\bm x)$, we set the problem
\begin{equation}\label{25}
 \operatorname{rot2\textbf{D}} \big (k(\bm x) \operatorname{rot2D} u \big ) + c(\bm x) u = f(\bm x), 
 \quad \bm x \in \Omega , 
\end{equation} 
\begin{equation}\label{26}
 u (\bm x) = 0,
 \quad \bm x \in \partial \Omega .
\end{equation} 
After approximating (\ref{25}), (\ref{26}) we come to the problem (\ref{19}) in which
in which the operator $A: H(\omega) \rightarrow H(\omega)$ is
\[
 A y = \operatorname{rot2\textbf{D}}_h \big (k(\bm x) \operatorname{rot2D}_h y \big ) + c(\bm x) y.
\]  
The self-adjointness property of this grid operator follows from the consistency of the approximations of the rotor operators in the sense of the fulfillment of the equality
\[
 (\operatorname{rot} \bm u, \bm v) - (\bm u, \operatorname{rot} \bm v) = 0
\] 
for functions (\ref{24}).

For two-dimensional problems (\ref{23}), (\ref{24}) with $\bm u(\bm x) = \{u_1, u_2\}$, the approximate solution at the boundary given (\ref{24}) is 
\begin{equation}\label{27}
 (\bm y \times \bm n) (\bm x) = 0,
 \quad \bm x \in \partial \omega^* . 
\end{equation} 
In the interior nodes of the $\omega^*$ grid, let's put
\[
 \operatorname{rot2D}_h \big (k(\bm x) \operatorname{rot2\textbf{D}}_h \bm y \big ) + c(\bm x) \bm y = \bm f(\bm x),
 \quad \bm x \in \omega^* - \partial \omega^*.
\]  
The grid problem is formulated in the form (\ref{22}) with operator $A: \bm H(\omega^*) \rightarrow \bm H(\omega^*)$, where
\[
 A \bm y = \operatorname{rot2D}_h \big (k(\bm x) \operatorname{rot2\textbf{D}}_h \bm y \big ) + c(\bm x) \bm y.
\] 
on the set of functions satisfying (\ref{27}). 

\section{Acknowledgment} 

The work was supported by the Russian Science Foundation (grant No. 24-11-00058).


\begin{thebibliography}{10}
\providecommand{\url}[1]{{#1}}
\providecommand{\urlprefix}{URL }
\expandafter\ifx\csname urlstyle\endcsname\relax
  \providecommand{\doi}[1]{DOI~\discretionary{}{}{}#1}\else
  \providecommand{\doi}{DOI~\discretionary{}{}{}\begingroup
  \urlstyle{rm}\Url}\fi

\bibitem{castillo2013mimetic}
Castillo, J.E., Miranda, G.F.: Mimetic Discretization Methods.
\newblock CRC Press (2013)

\bibitem{FVM}
Eymard, R., Gallou\"{e}t, T., Herbin, R.: Finite volume methods.
\newblock In: Handbook of Numerical Analysis, vol.~7, pp. 713--1020. North
  Holland, Amsterdam (2000)

\bibitem{frjazinov1975}
Fryazinov, I.V.: An approximation of mixed derivatives.
\newblock USSR Computational Mathematics and Mathematical Physics
  \textbf{15}(3), 95--111 (1975)

\bibitem{Generalized_difference_methods}
Li, R., Chen, Z., Wu, W.: Generalized Difference Methods for Differential
  Equations: Numerical Analysis of Finite Volume Methods.
\newblock Marcel Dekker (2000)

\bibitem{Samarskii1989}
Samarskii, A.A.: Theory of Difference Schemes.
\newblock Marcel Dekker, New York (2001)

\bibitem{samarskiiUMN}
Samarskii, A.A., Fryazinov, I.V.: Difference approximation methods for problems
  of mathematical physics.
\newblock Russian Mathematical Surveys \textbf{31}(6), 179--213 (1976)

\bibitem{SamarskiiKoldobaPoveschenkoFavorskii1996}
Samarskii, A.A., Koldoba, A.V., Povechenko, Y.A., Tishkin, V.F., Favorskii,
  A.P.: Difference Schemes on Unstructured Grids.
\newblock Kriterii, Minsk (1996).
\newblock In Russian

\bibitem{shashkov2018conservative}
Shashkov, M.: Conservative Finite-Difference Methods on General Grids.
\newblock CRC press (1996)

\bibitem{Strikwerda2004}
Strikwerda, J.C.: Finite Difference Schemes and Partial Differential Equations.
\newblock Society for Industrial Mathematics (2004)

\bibitem{TikhonovSamarskii1961}
Tikhonov, A.N., Samarskii, A.A.: Homogeneous difference schemes.
\newblock USSR Computational Mathematics and Mathematical Physics
  \textbf{1}(1), 5--67 (1962)

\bibitem{vabishchevich2005finite}
Vabishchevich, P.N.: Finite-difference approximation of mathematical physics
  problems on irregular grids.
\newblock Computational Methods in Applied Mathematics \textbf{5}(3), 294--330
  (2005)

\bibitem{vabishchevich2024gvm}
Vabishchevich, P.N.: Difference operator approximations on nonstandard
  rectangular grid.
\newblock Computational Mathematics and Mathematical Physics \textbf{64}(7),
  1367--1380 (2024)

\bibitem{da2014mimetic}
da~Veiga, L.B., Lipnikov, K., Manzini, G.: The Mimetic Finite Difference Method
  for Elliptic Problems.
\newblock Springer (2014)

\end{thebibliography}
\end{document}